\documentclass[preprint,12pt]{elsarticle}

\usepackage{amssymb}
\usepackage{amsmath}
\usepackage{subfigure}

\usepackage{placeins}

\usepackage{amsthm}
\newtheorem{remark}{Remark}
\newtheorem{definition}{Definition}

\DeclareMathOperator{\R}{\mathbb{R}}
\DeclareMathOperator{\XX}{\mathbb{X}}
\DeclareMathOperator{\TT}{\mathcal{T}}
\DeclareMathOperator{\pp}{\mathcal{P}}
\DeclareMathOperator{\Temp}{\mathfrak{T}}
\DeclareMathOperator{\Rey}{Re}
\DeclareMathOperator{\Ma}{Ma}
\DeclareMathOperator{\tol}{tol}

\newcommand{\ww}{\boldsymbol{w}}
\newcommand{\xx}{\boldsymbol{x}}

\journal{Computers \& Fluids}

\begin{document}

\begin{frontmatter}



\title{A temporally adaptive hybridized discontinuous Galerkin method for time-dependent compressible flows}



\author[jochen]{Alexander Jaust}
\author[jochen]{Jochen Sch\"utz\corref{tel}}

\address[jochen]{Institut f\"ur Geometrie und Praktische Mathematik, RWTH Aachen University, Aachen, Germany}
\cortext[tel]{Phone: (+49) 241 80 97677, Fax: (+49) 241 80 92349}
\ead{schuetz@igpm.rwth-aachen.de}


\begin{abstract}
 The potential of the hybridized discontinuous Galerkin (HDG) method has been recognized for the computation of stationary flows. Extending the method to time-dependent problems can, e.g., be done by backward difference formulae (BDF) or diagonally implicit Runge-Kutta (DIRK) methods. In this work, we investigate the use of \emph{embedded} DIRK methods in an HDG solver, including the use of adaptive time-step control. 
 Numerical results demonstrate the performance of the method for both linear and nonlinear (systems of) time-dependent convection-diffusion equations.

\end{abstract}

\begin{keyword}
hybridized discontinuous Galerkin method \sep embedded diagonally implicit Runge-Kutta methods \sep time-dependent convection-diffusion equations


\end{keyword}

\end{frontmatter}

\section{Introduction}\label{sec:introduction}
The last few years have seen a tremendous increase in the use and development of high order methods for aerodynamic applications, see \cite{Wang2003,Huynh09,Hartmann2002a} to only mention a few. These methods represent the unknown function by a (piecewise) polynomial of degree larger than two, and so exceed the design order of Finite-Volume schemes that are nowadays standard tools in the aerospace industry. One particular example of high order methods is the discontinuous Galerkin (DG) method introduced by Reed and Hill \cite{ReHi73} and subsequently extended to all sorts of equations by many authors, see, e.g., \cite{Bassi1997,DG1,DG2,DG3,DG4,DG5,BaOd00,ABCM}. DG offers a lot of inherent advantages, such as the flexibility concerning the local degree of polynomial, the easy incorporation of boundary conditions on complicated domains, local conservativity and many more. 

However, in the context of stationary problems, where usually the Jacobian of the method is needed, DG methods suffer from large memory requirements. This is due to the fact that the number of unknowns increases as $O((p+1)^{d} N)$, where $p$ is the order of the local polynomial, $d$ the spatial dimension and $N$ the number of elements in the triangulation. Especially for the combination of high $p$ and $d > 1$, this poses a severe restriction. One way of tackling this problem is to use \emph{hybridized} DG (HDG) methods, see \cite{NPC09L,NPC09,PeNgCo,NgPeCo11,NgPe12,COGOLA,EgSch09,SchMa11,WoMaSc13}. The globally coupled unknowns in this type of method are the degrees of freedom belonging to the unknown function $w$ on the \emph{edges} instead of on the elements. Obviously, this yields a reduction in dimension, and the number of globally coupled unknowns behaves as $O((p+1)^{d-1} \widehat N)$, where $\widehat N$ is the number of edges in the triangulation. 

During the second high-order workshop at DLR in Cologne (see \cite{HighOrderCFD} for a summary of the first workshop), we have presented our method for 'easy' stationary problems in the context of two-dimensional Euler and Navier-Stokes equations. It could be seen that hybridized DG methods have a potential of outperforming more traditional schemes. 
Still, there are a few things missing, among them the efficient implementation of time-integration routines. At least conceptually, the temporal discretization for DG schemes can be done using a method of lines approach. Unfortunately, this is not possible for the hybrid method. However, using a dual time-stepping approach \cite{Jameson1991}, one can incorporate \emph{implicit} time integration methods. 


Extension to time-dependent problems has been made using BDF (backward differentiation formulae)  methods \cite{NPC09L,NPC09,SchWoMay2012} and DIRK (diagonally implicit Runge-Kutta) schemes \cite{NgPeCo11, NgPe12}. 
In this work, we compare different DIRK schemes \cite{HaiWan, Rabeh87,  Cash} and investigate their use for the temporal discretization of the hybridized DG method including time-step control. As the HDG method applied to a temporal problem gives rise to a differential algebraic equation \cite{HaiWan}, the DIRK schemes have to be chosen suitably. 
The goal of this paper is to demonstrate that the combination of two well-known ingredients (embedded DIRK and HDG) is possible in the context of compressible fluid flows. In this sense, this work constitutes an intermediate step that an efficient solver for unsteady aerodynamic problems can also be based on HDG. 
We show numerical results for advection-diffusion, Euler and Navier-Stokes equations in two dimensions, demonstrating that the combination of HDG with embedded DIRK methods yields a stable method. 

The paper is organized as follows: In Section \ref{sec:equations}, we shortly introduce the underlying equations. In Section \ref{sec:method}, we introduce the HDG method and discretize the resulting semi-discrete system using embedded DIRK methods. In addition, time-step control is discussed in this section. Section \ref{sec:numerical_results} shows numerical results, and Section \ref{sec:conclusions} offers conclusions and an outlook. In the appendix section we give the Butcher tableaus of the embedded DIRK methods we use. 

\section{Underlying Equations}\label{sec:equations}
In this work, we consider for a two-dimensional domain $\Omega \subset \R^2$ the general unsteady convection-diffusion equation, given as 
\begin{alignat}{2}
 \label{eq:conv_diff}
 w_t + \nabla \cdot \left(f(w) - f_v(w, \nabla w)\right) &= g &\quad& \forall (x, t) \in \Omega \times (0, T),  \\
 w(x, 0) &= w_0(x) &\quad& \forall x \in \Omega,
\end{alignat}
where $w_0$ and $g$ are given functions, and $f : \R^m \rightarrow \R^{m \times 2}$, $f_v : \R^m \times \R^{m \times 2} \rightarrow \R^{m \times 2}$ are convective and diffusive flux, respectively. $T > 0$ denotes a final time; while $m$ is the dimension of the system. 
The equations are equipped with appropriate boundary conditions, which depend on the particular choice of the fluxes $f$ and $f_v$. 

Note that both unsteady Euler and Navier-Stokes equations fall into this framework, where the unknown is $w = (\rho, \rho u, \rho v, E)$, i.e., density, momentum, and total energy. Corresponding fluxes are defined by 
\begin{alignat}{3}
 \label{eulerfluxes}
 f_1 &= (\rho u, p + \rho u^2, \rho u v, u(E+p))^T, &\quad&
 f_2 &\ =\ & (\rho v, \rho u v, p + \rho v^2, v(E+p))^T,  \\
  f_{v}^1 &= (0, \tau_{11}, \tau_{21}, \tau_{11}u + \tau_{12}v + k \Temp_{x_1})^T, &\quad&
  f_{v}^2 &\ =\ & (0, \tau_{12}, \tau_{22}, \tau_{21}u + \tau_{22}v + k \Temp_{x_2})^T, 
\end{alignat}
and the right-hand side $g \equiv 0$. For Euler equations, $f_v \equiv 0$. As usual, $p$ denotes pressure, $\tau$ stress tensor, $\Temp$ temperature and $k$ thermal conductivity coefficient.
$p$ is coupled to the conservative variables using the ideal gas law in the form 
\begin{align}
 p = (\gamma - 1) \left(E - \frac{\rho (u^2 + v^2)}{2}\right). 
\end{align}
For a wide range of flow conditions, the ratio of specific heats $\gamma$ is assumed to be constant with value 1.4. If not stated otherwise, we rely on this choice. 

As it is frequently done when considering diffusion equations \cite{ABCM}, we formulate \eqref{eq:conv_diff} as a first-order system by introducing the unknown function $\sigma := \nabla w$, i.e., in the sequel, we consider
 \begin{alignat}{5}
 \label{eq:conv_diff_sigma}
 \sigma &= \nabla w, &\quad& & w_t + \nabla \cdot \left( f(w) - f_v(w, \sigma) \right) &=& g \hspace{.8cm} &\quad& \forall (x, t) \in \Omega \times (0, T), \\
  & & & & \hspace{3.5cm} w(x,0) &=&\ w_0(x)  &\quad& \forall x \in \Omega. 
\end{alignat}

\section{A Hybridized DG Method}\label{sec:method}
\subsection{Semi-Discrete Method}
In this section, we discretize equation \eqref{eq:conv_diff} in space using the hybridized DG method.
To this end, we need a triangulation that is defined in the sequel:

\begin{definition}
We assume that $\Omega$ is triangulated as 
\begin{align}
 \Omega = \bigcup_{k=1}^N \Omega_k. 
\end{align}
We define an edge $e_k$ to be either an intersection of two neighboring elements, or the intersection of an element with the physical boundary $\partial \Omega$, having positive one-dimensional measure. $\Gamma$ denotes the collection of all these intersections, while $\Gamma_0 \subset \Gamma$ denotes those $e_k \in \Gamma$ that do not intersect the physical boundary $\partial \Omega$ of the domain. We define $\widehat N := |\Gamma|$ to be the number of edges in $\Gamma$.   
\end{definition}

For the ease of presentation, we introduce the following standard abbreviations for integration: 
\begin{alignat}{3}
 (f_1,f_2) &:= \sum_{k=1}^N \int_{\Omega_k} f_1 \cdot f_2 \ dx, &\quad&
 \langle f_1,f_2 \rangle_{\Gamma} := \sum_{k=1}^{\widehat N} \int_{e_k} f_1 \cdot f_2 \ d \sigma, &\quad& \\ \langle f_1,f_2 \rangle_{\partial \Omega_k} &:= \sum_{k=1}^{N} \int_{\partial \Omega_k} f_1 \cdot f_2 \ d \sigma. 
\end{alignat}
In the method to be presented, both $\sigma$ and $w$ are approximated explicitly. Additionally, we introduce a variable $\lambda$ that has support on the skeleton of the mesh only, $\lambda := w_{|\Gamma_0}$. The resulting algorithm will thus approximate the quantity 
\begin{align}
  \label{eq:ww}
  \ww := (\sigma, w, w_{|\Gamma_0}). 
\end{align}
On first sight, this seems like a tremendous increase in degrees of freedom, as one does not only approximate $w$ (which is usually done in DG methods), but also $\sigma$ and $\lambda$. However, the hybridized DG algorithm is constructed in such a way that one can locally eliminate both approximations to $\sigma$ and $w$ in favor of the approximation to $\lambda$, see \cite{CoGo04}. 
The only coupled degrees of freedom are those associated to the approximation of $\lambda$. 

In the sequel, we define the correct approximation spaces: 
\begin{definition}
 Let the approximation to $\ww(\cdot, t)$ at some fixed time $t$,
 \begin{align}
    \label{eq:ww_h}
    \ww_h(\cdot, t) := (\sigma_h(\cdot, t), w_h(\cdot, t), \lambda_h(\cdot, t)) 
 \end{align}
 be in $\XX_h := H_h \times V_h \times M_h$, where 
 \begin{alignat}{4}
  H_h &:= \{ f \in L^2(\Omega) \ &|& \ f_{|\Omega_k} \ &\in& \ \Pi^p(\Omega_k) \quad &&\forall k = 1, \ldots N\}^{2 \cdot m} \\
  V_h &:= \{ f \in L^2(\Omega) \ &|& \ f_{|\Omega_k} \ &\in& \ \Pi^p(\Omega_k) \quad &&\forall k = 1, \ldots N\}^m \\
  M_h &:= \{ f \in L^2(\Gamma) \ &|& \ f_{|e_k} \ &\in& \ \Pi^p(e_k) \quad &&\forall k = 1, \ldots \widehat N, e_k \in \Gamma\}^m.
\end{alignat}
\end{definition}
\begin{remark}
 Whenever we use bold letters for a function, we think of a triple of functions. See the definitions \eqref{eq:ww} and \eqref{eq:ww_h} of $\ww$  and $\ww_h$, respectively, for examples.
\end{remark}

Based on these approximation spaces, and following Nguyen et al.'s and our previous work \cite{NPC09L, NPC09, SchMa11}, we can define the semi-discretization in the sequel:
\begin{definition}
 A semi-discrete approximation 
 \begin{align}
  \ww_h(\cdot, t) := (\sigma_h(\cdot, t), w_h(\cdot, t), \lambda_h(\cdot, t)) \in \XX_h
 \end{align}
 to \eqref{eq:conv_diff_sigma} using the hybridized DG method is defined as the function $\ww_h$, such that for all $t \in (0, T)$: 
 
\small
\begin{alignat}{2}  
  \label{eq:hybrid1}
  (\sigma_h - \nabla w_h, \tau_h) - \langle \lambda_h - w_h^-, \tau_h^- \cdot n \rangle_{\partial \Omega_k} &= 0  &\quad& \forall \tau_h \in H_h\\
  \label{eq:hybrid2}
  ((w_h)_t, \varphi_h) - (f(w_h) - f_v(w_h, \sigma_h), \nabla \varphi_h) \\ + \langle (\widehat f - \widehat f_v)  \cdot n, \varphi_h^- \rangle_{\partial \Omega_k} &= (g, \varphi_h)  &\quad& \forall \varphi_h \in V_h\\
  \label{eq:hybrid3}
  \langle  [(\widehat{f} - \widehat{f_v})] \cdot n , \mu_h \rangle_{\Gamma} &= 0 &\quad& \forall \mu_h \in M_h.
\end{alignat}
\normalsize
Numerical fluxes $\widehat f$ and $\widehat f_v$ are defined as 
\begin{alignat}{2}
 \widehat {f}   &:= f(\lambda_h) - \delta \left(\lambda_h - w_h^- \right) n, &\quad
 \widehat {f_v} &:= f_v(\lambda_h, \sigma_h^-) + \tau \left(\lambda_h - w_h^- \right) n.
\end{alignat}
Both $\delta$ and $\tau$ are real parameters that depend on $f$ and $f_v$, respectively. 
\end{definition}
\begin{remark}
In the limiting cases $f \equiv 0$, one has $\delta = 0$, while for $f_v \equiv 0$, $\tau = 0$. Boundary conditions are incorporated into the definition of the fluxes $\widehat f$ and $\widehat{f_v}$, i.e., the definition of $\widehat f$ and $\widehat f_v$ is altered on $\Gamma \backslash \Gamma_0$. This is done in such a way that the method is \emph{adjoint consistent}. For the ease of presentation, we neglect the details which can be found in, e.g., \cite{SchMa11, SM12}. 
The fluxes are such that $w_h^+$ and $w_h^-$ (and $\sigma_h^+$ and $\sigma_h^-$, respectively), are not directly coupled. 
This allows for a static condensation, so that the only globally coupled degrees of freedom are those associated to $\lambda_h$. This is in general a reduction of degrees of freedom in comparison to traditional DG methods. 
\end{remark}

\begin{remark}
 In the stationary case, numerical results show \cite{SchMa11AIAA} that the expected optimal convergence rates $p+1$ are met for both the $L^2-$error in $w_h$ and $\sigma_h$. 
\end{remark}

For the ease of presentation, we rewrite \eqref{eq:hybrid1}-\eqref{eq:hybrid3} as 
\begin{alignat}{2}
 \label{eq:hybrid}
 \TT((w_h)_t, \varphi_h)  + N(\ww_h; \xx_h) = 0 &\quad& \forall \xx_h \in \XX_h,
\end{alignat}
where $\xx_h := (\tau_h, \varphi_h, \mu_h)$ is just a convenient shortcut for the test functions. 
$\TT$ denotes the vector having $0$ entries for the first and the last equation, i.e., 
\begin{align}
 \TT((w_h)_t, \varphi_h) :=  (0, ((w_h)_t, \varphi_h), 0)^T, 
\end{align}
and $N(\ww_h; \xx_h)$ is the remaining part belonging to the discretization of the stationary convection-diffusion equation.
\begin{remark}
 A straightforward method of lines approach can only be applied if $\TT$ does not have the zero entries in the first and the last argument. At least in principle, one could derive equations for $\sigma_t$ and $\lambda_t$ and try to incorporate them. However, this would require a large amount of derivatives which will most likely deteriorate the order of the scheme. To this end, we will in the sequel rely on \emph{implicit} time discretization methods and Jameson's idea of dual time-stepping \cite{Jameson1991}. 
\end{remark}

\subsection{(Embedded) DIRK Discretization}
In this section, we explain how time is discretized in our setting, so that in the end we get a fully discrete algorithm. In an earlier work \cite{SchWoMay2012}, we have used backward difference formulae for the discretization of the temporal part. These methods seemed particularly suited to the method at hand and showed very good accuracy results. 
However, they suffer from the need to use expensive initial step(s) to obtain suitable start values and the complicated incorporation of temporal adaptivity. Furthermore, the extension to higher order is difficult, because those methods cease to be A-stable for order of consistency larger than 2, which can actually cause stability problems. In the current paper, inspired by the recent work of Nguyen, Peraire and Cockburn \cite{NgPeCo11} and Nguyen and Peraire \cite{NgPe12}, we use diagonally implicit Runge-Kutta (DIRK) methods and, more specifically, \emph{embedded} DIRK methods to achieve an adaptive temporal discretization of (potentially) high order. 

We consider an adaptive sequence of time instances $0 = t^0 < t^1 < \ldots < t^M = T$, where both $\Delta t^n := t^{n+1}-t^n$ and $M$ depend on the solution to be approximated. We define $\ww_h^n$ to be the approximation of $\ww_h(\cdot, t^n)$ at time instance $t^n$.

As \eqref{eq:hybrid} constitutes a differential algebraic equation (DAE), special care has to be taken in choosing suitable methods. Desirable properties are A-stability to allow for large time steps; and the DIRK schemes should have no explicit stage, which is a necessary condition for convergence in the case of DAEs, see also Remark \ref{rem:sdirk}. 
A-stable methods that fulfill this condition can be rarely found in the literature, we use those presented in \cite{HaiWan, Rabeh87, Cash}. 

\begin{definition}
 An embedded DIRK method is given by its Butcher tableau with a lower triangular matrix $A \in \R^{k \times k}$, a node vector $\beta \in \R^k$ and two weighting vectors $\gamma_1, \gamma_2 \in \R^k$. Frequently, the Butcher tableau is given as 
\begin{align}
 \begin{array}{l|llll}
  \beta_1 & A_{11} \\ 
  \beta_2 & A_{21} & A_{22} \\
  \vdots  & \vdots & \vdots &  \ddots \\
  \beta_k & A_{k1} & \hdots & \hdots & A_{kk} \\ \hline
    & \gamma_{11} &\ldots &\ldots &\gamma_{1k} \\ \hline
    & \gamma_{21} &\ldots &\ldots &\gamma_{2k} \\ \hline
 \end{array}.
\end{align}
\end{definition}
\begin{remark}\label{rem:dirk_ode}
 Applying an embedded DIRK method to an ordinary differential equation 
 \begin{align}
  y'(t) = f(t, y(t))
 \end{align}
 amounts to approximating two values $y_1^{n+1}$ and $y_2^{n+1}$. Those two values are given by 
 \begin{alignat}{2}
  \label{eq:dirk1}
  y_1^{n+1} &= y^n + \Delta t^n \sum_{i=1}^k \gamma_{1i} f(t^n + \beta_i \Delta t^n, y^{n, i}) \\
  \label{eq:dirk2}
  y_2^{n+1} &= y^n + \Delta t^n \sum_{i=1}^k \gamma_{2i} f(t^n + \beta_i \Delta t^n, y^{n, i})
 \end{alignat}
 with the intermediate values $y^{n, i}$ implicitly given by 
 \begin{align}
  y^{n, i} = y^n + \Delta t^n \sum_{j=1}^i A_{ij} f(t^n + \beta_j \Delta t^n, y^{n, j}).
 \end{align}
 Note that each step of computing the $y^{n, i}$ is basically an implicit Euler step. 
\end{remark}
\begin{remark}\label{rem:embedded}
 Note that the two approximations to $y(t^{n+1})$ are supposed to have different orders of accuracy. More precisely, we have chosen to enumerate in such a way that $y_1^{n+1}$ is the more accurate approximation. This means that $\|y_1^{n+1} - y_2^{n+1}\|$ can serve as a measure of consistency error. 
 Furthermore, we choose our schemes in such a way that the DIRK method corresponding to $\gamma_1$ is A-stable. Usually, $\gamma_{2k}$ is zero, so that the embedded method has actually only $k-1$ stages. 
\end{remark}


It is straightforward to apply the DIRK method to fully discretize \eqref{eq:hybrid}: As in Rem. \ref{rem:dirk_ode}, $\ww_h(\cdot, t^{n+1})$ is approximated by two different values $\ww_{h,1}^{n+1}$ and $\ww_{h, 2}^{n+1}$. The approximations are such that 
\begin{alignat}{2}
 \label{eq:hybriddirk1}
 \TT\left(w_{h, 1}^{n+1}-w_h^n, \varphi_h\right) + \Delta t^n \sum_{i = 1}^k \gamma_{1i} N(\ww_h^{n,i}; \xx_h) &= 0 &\quad& \forall \xx_h \in \XX_h\\
 \label{eq:hybriddirk2}
 \TT\left(w_{h, 2}^{n+1}-w_h^n, \varphi_h\right) + \Delta t^n \sum_{i = 1}^k \gamma_{2i} N(\ww_h^{n,i}; \xx_h) &= 0 &\quad& \forall \xx_h \in \XX_h. 
\end{alignat}
The intermediate stages $\ww_h^{n,i}$ are defined via the equation 
\begin{alignat}{2}
 \label{eq:hybriddirk_intermediate}
 \TT\left({w_h^{n, i} - w_h^n}, \varphi_h\right) + {\Delta t^n}\sum_{j=1}^i A_{ij} N(\ww_h^{n, j}; \xx_h) &= 0 &\quad& \forall \xx_h \in \XX_h. 
\end{alignat}
\begin{remark}
 Note that the computation of $\ww_h^{n, i}$, see \eqref{eq:hybriddirk_intermediate}, reduces to solving a stationary system, and thus fits very nicely into the framework of our steady-state solver \cite{SchMa11}. Note furthermore that \eqref{eq:hybriddirk1}-\eqref{eq:hybriddirk2} is an explicit step (up to the inversion of a mass matrix) and corresponds to \eqref{eq:dirk1}-\eqref{eq:dirk2}.
\end{remark}

\begin{remark}
 In our numerical computations, we use three different embedded DIRK schemes, one taken from the book by Hairer and Wanner \cite{HaiWan}, one taken from Al-Rabeh \cite{Rabeh87} and one taken from Cash \cite{Cash}. The corresponding tableaus are given in \ref{sec:appendix_dirk}. The design orders are $q = 3$ for Cash's method and $q = 4$ for the other two methods. 
\end{remark}

\begin{remark}\label{rem:sdirk}
 All the embedded DIRK schemes fall within the class of SDIRK (i.e., singly DIRK) which in particular means that all the diagonal entries of the Butcher tableau are different from zero. As the semi-discrete system \eqref{eq:hybrid} is a differential-algebraic equation, this is necessary to ensure stability \cite{HaiWan}. 
\end{remark}

As we are using an embedded DIRK scheme, it is our desire to adaptively control the time-step $\Delta t^n$. To this end, we define the following error estimation based on the quantities $w_{h,1}^{n+1}$ and $w_{h, 2}^{n+1}$:
\begin{align}
e_h^n := \|w_{h,1}^{n+1}-w_{h,2}^{n+1}\|_{L^2}.
\end{align}
Note that the use of a non-bold $w$ is not a typo, we only use the second component of $\ww_{h, i}^{n+1}$ which represents the solution within the elements, and is probably the most important quantity.
As is customary in the use of embedded Runge-Kutta methods, a time-step is rejected (i.e., repeated with a fraction of the time-step size, e.g., half the time-step) if 
\begin{align}
 e_h^n > \Delta t^n \cdot \tol 
\end{align}
for a user-defined tolerance $\tol$. This approach guarantees that $\sum_{i = 0}^{M-1} e_h^n < T \cdot \tol$. 
If the time-step is accepted, we take $\ww_{h,1}^{n+1}$ to be the new approximate value $\ww_h^{n+1}$ and compute the new time-step, based on the old time-step, as
\begin{alignat}{2}
  \Delta t^{n+1} = \alpha \Delta t^n \left(r_h^n\right)^{-\frac{1}{q-1}}.
\end{alignat}
$\alpha$ is a safety factor given by 
\begin{alignat}{2}
 \alpha = 0.9 \frac{2 n_{it,max} + 1}{2 n_{it,max} + n_{it}}, &\quad&
\end{alignat}
where we take the maximum number of Newton steps per stage $n_{it}$, and the maximum allowable number of Newton steps $n_{it, max}$ in our nonlinear solver into account. $r_h^n$ is defined by
\begin{align}
 r_h^n := \frac{e_h^n}{\tol \Delta t^n}. 
\end{align}
The underlying paradigm is that $r_h^n$ is close to one, because if $e_h^{n+1} \approx e_h^n$, this allows the largest possible time-step such that $e_h^{n+1}$ is not too big. 

For a detailed derivation, see standard textbooks such as \cite{HaiWan}. The design accuracy of the  DIRK scheme is denoted by $q$. As usual, $\Delta t^n$ is 'limited' such that it does not exceed a maximum and a minimum value.  

\begin{remark}
 We control the higher-order Runge-Kutta method with the lower-order one. Strictly speaking, it should be the other way around. Nevertheless, it has been done in literature (see, e.g., \cite{HaiWan}) as it is desirable to keep the higher-order approximation, and we do it for the same reason. 
\end{remark}


\section{Numerical Results}\label{sec:numerical_results}
In this section, we show numerical results obtained with our method. We consider two-dimensional problems. First, we start from a simple, convection-diffusion problem to test the accuracy of our method and also the performance of the time-step control. Then, we show results for both Euler and Navier-Stokes equations. 

In these numerical examples, the nonlinear system of equations \eqref{eq:hybriddirk_intermediate} is solved using a damped Newton's method, thereby obtaining a sequence of linear systems of equations in the variations of $\ww_h^{n, j}$.
Then, one can apply static condensation such that one only has to solve for the variations in $\lambda_h^{n,j}$  (see \cite[III.C.]{SchMa11AIAA} for a more detailed derivation of Newton's method in this context). This resulting linear system is solved using an ILU(0) preconditioned GMRES through the PETSc \cite{petsc1,petsc2,petsc3} library with a relative tolerance of $10^{-4}$. 
Obtaining the variations in $w_h^{n,j}$ and $\sigma_h^{n,j}$ necessitates the solution of multiple (comparably) small linear systems of equations. These systems are solved with LAPACK routine dgesv \cite{lapack}. Convergence of a nonlinear iteration is obtained if the 2-norm of the residual of the equation associated with $\lambda_h^{n,j}$ drops below $10^{-10}$. 

\subsection{Scalar Convection-Diffusion Equation (Rotating Gaussian)}

As a scalar convection-diffusion equation, we present a test case that has previously been investigated by Nguyen et al. \cite{NPC09L}. The problem is both scalar and linear, with convective and viscous flux vector, respectively, given as
\begin{alignat*}{2}
 f(w) &= (-4y, 4x)^T w, &\quad& f_v(w, \nabla w) = 10^{-3} \nabla w.
\end{alignat*}
The source term $g$ is set to zero, and we consider the domain $\Omega = [-0.5,0.5]^2$. The final time $T$ is defined as $T = \frac\pi4$. Note that this test case is very interesting as, in the vicinity of the origin, the problem is diffusion-dominated, while, away from the origin, it is convection-dominated. Initial data are given by a (scaled) Gaussian distribution, i.e.,
\begin{align}
 w(x, y, 0) := e^{-50 (x^2 + y^2)}. 
\end{align}
An exact solution to this problem is known, and on $\partial \Omega$, we impose Dirichlet boundary conditions that we choose to be this exact solution. In the numerical results, we use the $L^2-$ norm of $w-w_h$ at final time $T$ as a measure of error. 

\begin{figure}[htb]
 \begin{center}
  \includegraphics{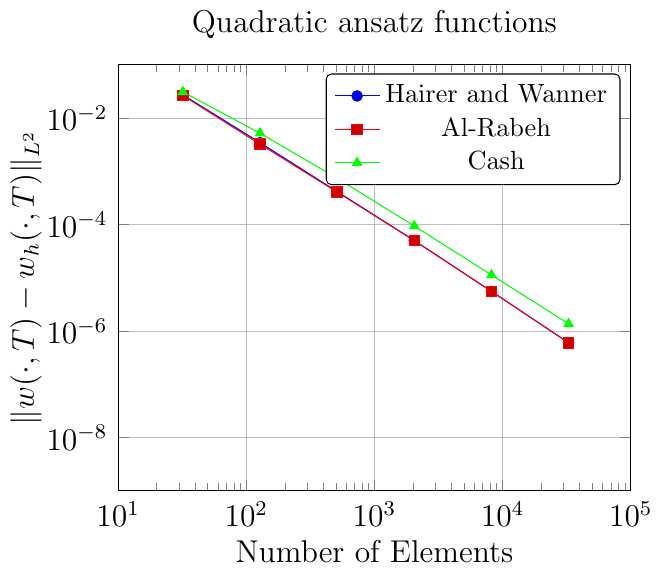}
  \includegraphics{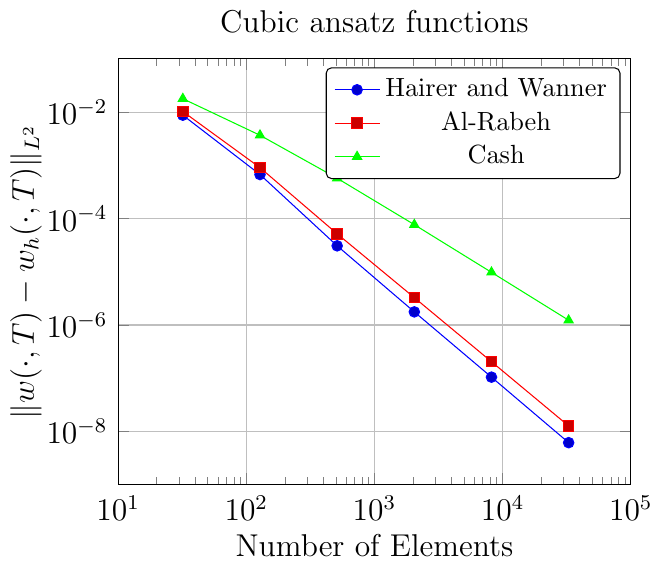}
  \caption{Rotating Gaussian: Convergence of $\|w(\cdot, T)-w_h(\cdot, T)\|_{L^2}$ under uniform spatial and temporal refinement.}	
  \label{fig:convg_linear_uniform}
 \end{center}

\end{figure}

In Fig. \ref{fig:convg_linear_uniform}, we demonstrate that the fully discrete scheme (without time-step control) converges under both uniform spatial and temporal refinement. The design accuracy of Cash's DIRK scheme is $q = 3$, while the design accuracy of the other two DIRK schemes is $q = 4$. The convergence order to be expected is thus $\min(q, p+1)$. For piecewise quadratic polynomials as ansatz functions, i.e., for $p = 2$, one can thus expect third order of convergence, while for piecewise cubic polynomials, i.e., for $p = 3$, one can expect third and fourth order of convergence, respectively. Numerical results confirm this expectation. 

\begin{figure}[h]
 \begin{center}
  \subfigure[Uniform temporal refinement.]{\includegraphics{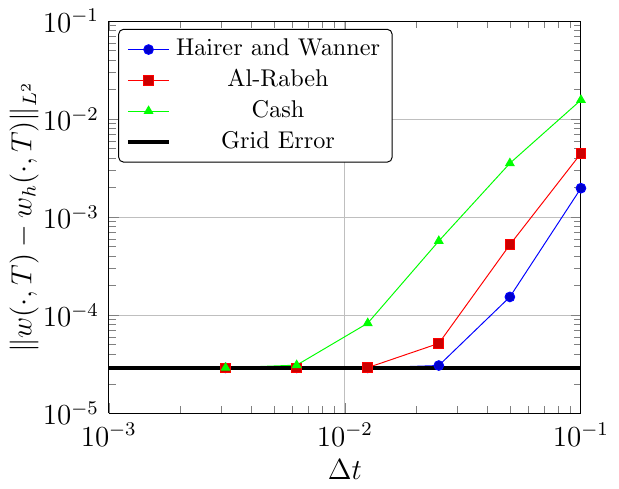}\label{fig:dt_cvg}}
  \subfigure[Adaptive temporal refinement.]{\includegraphics{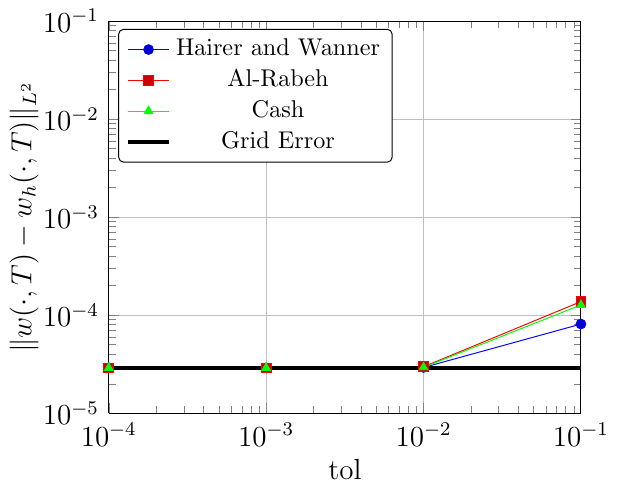}\label{fig:tol_cvg}}
  \caption{Rotating Gaussian: Convergence of $\|w(\cdot, T)-w_h(\cdot, T)\|_{L^2}$ for a fixed mesh with $N = 512$ elements and cubic ansatz functions.}
 \end{center}
\end{figure}

\begin{figure}
 \begin{center}
  \subfigure[$\tol = 10^{-1}$.]{\includegraphics{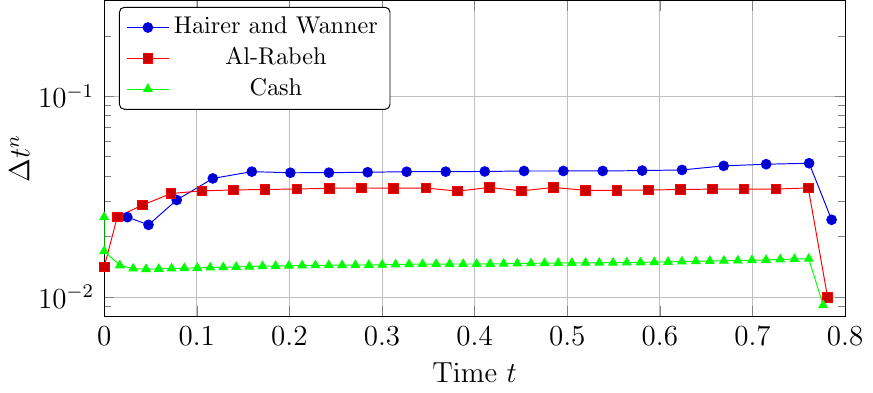}\label{fig:timestep_tol_1e1}}
  \subfigure[$\tol = 10^{-2}$.]{\includegraphics{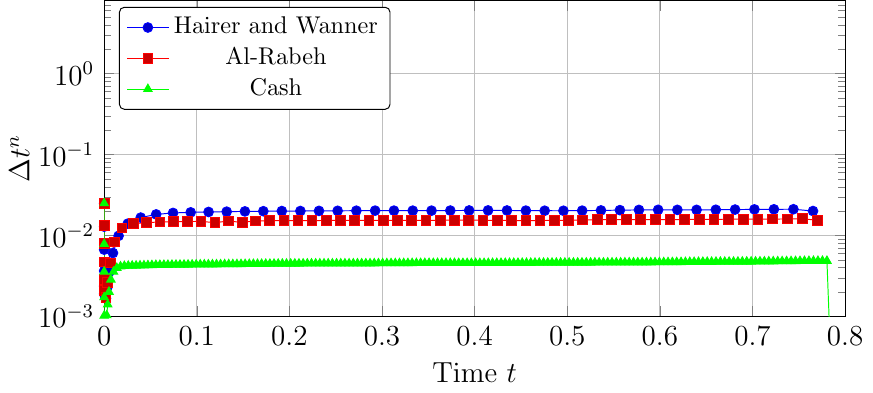}\label{fig:timestep_tol_1e2}}
  \caption{Rotating Gaussian: Time-step size for a fixed mesh with $N = 512$ elements and cubic ansatz functions for different values of $\tol$.}
 \end{center}

\end{figure}

The next test is about the time-step adaptation per se. We take a fixed spatial mesh consisting of 512 elements and cubic ansatz functions, and only refine in time. For $\Delta t \rightarrow 0$, this will yield the spatial error. In Fig. \ref{fig:dt_cvg}, we plot the error evolution for a fixed time-step $\Delta t$, while in Fig. \ref{fig:tol_cvg}, we plot the error versus different tolerances for all three DIRK methods. One can clearly see that all the methods are able to obtain the spatial error with only a moderate degree of tolerance. The quasi-optimal time-step is determined automatically. Furthermore, it can be observed that the method by Hairer and Wanner has the best convergence properties, both uniform and adaptive. In Figs. \ref{fig:timestep_tol_1e1}-\ref{fig:timestep_tol_1e2}, we plot the evolution of the time-step for $\tol = 10^{-1}$ and $\tol = 10^{-2}$. 
The test case under consideration is actually quite homogeneous in its temporal behavior. As a consequence, one can see that after an initial increase, $\Delta t^n$ remains nearly constant. 
Also here, Hairer and Wanner's method performs best, as it yields the largest time-step without sacrificing accuracy. 
The kink which can can be observed in the plots of the time-step size at the right is not an artifact, but in fact the result of a fixed final time $T$ which has to be reached by the last time-step.

\begin{figure}
 \begin{center}
  \includegraphics{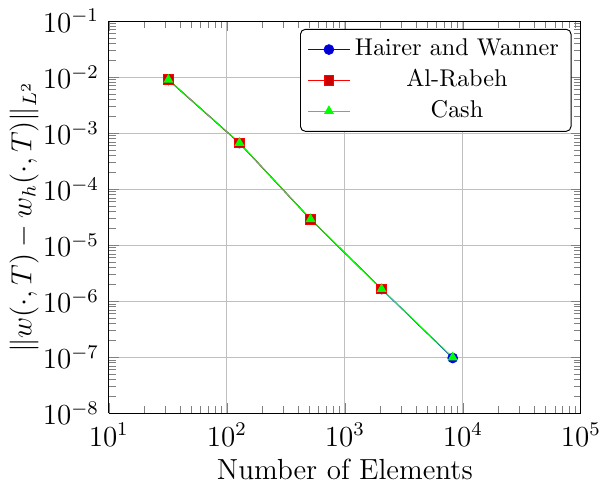}
 \end{center}
 \caption{Rotating Gaussian: Convergence of $\|w(\cdot, T)-w_h(\cdot, T)\|_{L^2}$ for cubic polynomials under both mesh and tolerance refinement.}\label{fig:convg_adaptive}
\end{figure}

Using a fixed tolerance $\tol$ is obviously not enough if one considers mesh refinement. With an increasing spatial resolution, the tolerance should decrease. We choose to set $\tol = O(h^{\min(q,p+1)})$, where $h \propto \frac{1}{\sqrt{N_x}}$ is a measure of the mesh size. We start on an initial mesh with $N = 32$ elements and an initial tolerance $\tol = 10^{-1}$. Convergence results can be seen in Fig. \ref{fig:convg_adaptive}. The surprising outcome of this is that the error values nearly lie on top of each other. This means that the time-step adaptation really performs well, as it obviously minimizes the temporal error. This also explains why results associated to Cash's method (which is only third order accurate) show fourth order behavior: With the chosen value of $\tol$, temporal error is not dominant any more, and all one observes is the influence of the spatial error.
Comparing with Fig. \ref{fig:convg_linear_uniform}, one can see that this adaptive approach performs better than even the uniform refinement. 

\subsection{Euler equations (Radial Expansion Wave)}

\begin{figure}
 \begin{center}
  \includegraphics[width=0.48\textwidth]{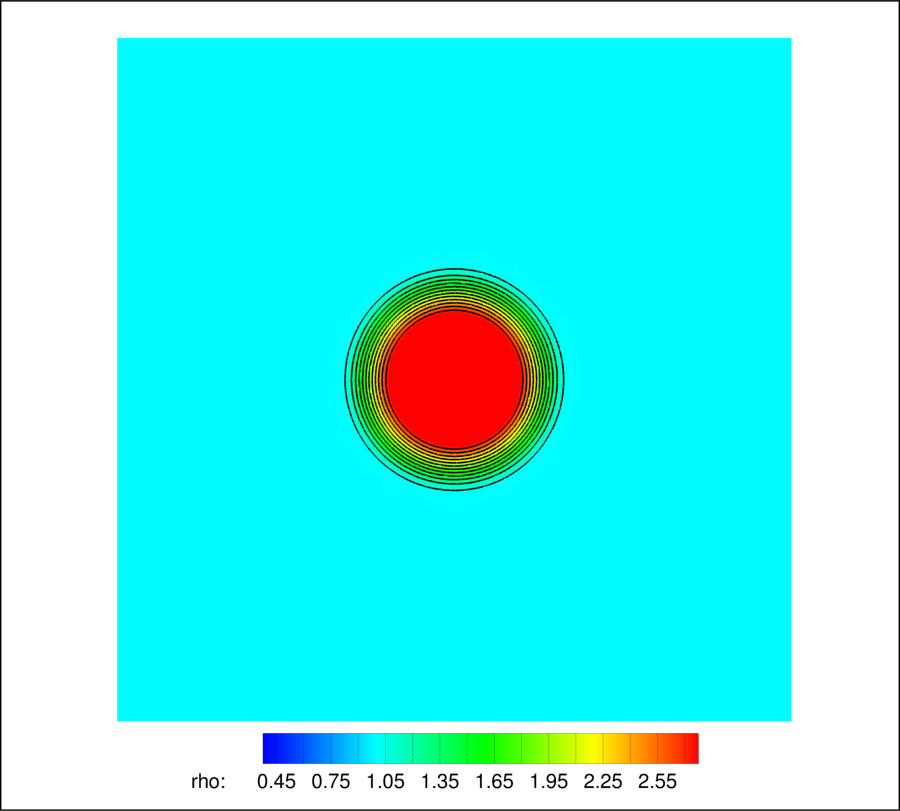}
  \includegraphics[width=0.48\textwidth]{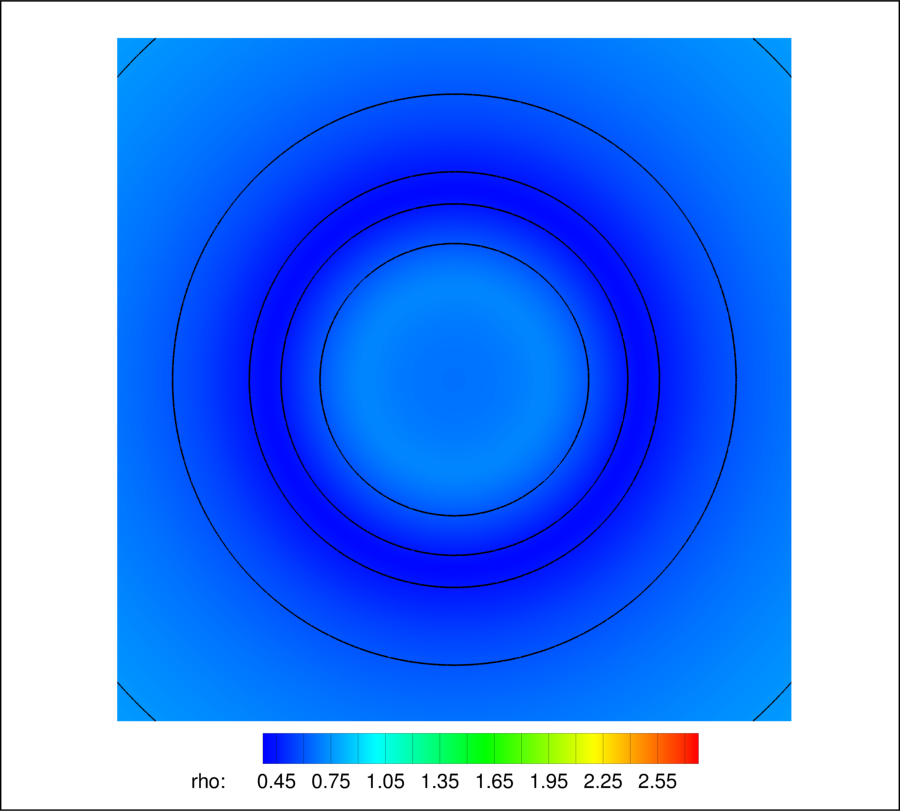}
 \end{center}
 \caption{Density of radial expansion wave at time $t = 0$ and $t = 2$, respectively.}\label{fig:radexpwave}
\end{figure}

The next test case has been proposed for both first and second high order workshop \cite{HighOrderCFD}. It is to compute a radially symmetric, inviscid flow using the Euler equations on domain $\Omega = [-4,4]^2$. For the standard choice of $\gamma = 1.4$, the flow ceases to be smooth and its derivative becomes discontinuous. To this end, it has been proposed to use $\gamma = 3$, which will be our choice in the sequel. The final time $T$ is set to $T = 2$.  The flow is supersonic throughout the domain, so it is fully specified by its initial conditions (for simplicity, $r \equiv r(x,y) := \sqrt{x^2 + y^2}$ denotes radius):
\begin{alignat}{4}
 q(r) &:= \begin{cases} 0, & 0 \leq r < \frac12 \\ 
		      \frac1{\gamma} \left(1 + \tanh\left(\frac{r-1}{0.25-(r-1)^2} \right) \right), & \frac12 \leq r < \frac32 \\
		      \frac2{\gamma}, &r \geq \frac32\end{cases}, \\
 u(x,y,0) &:= \frac{x}{r} q(r),  \quad \quad
 v(x,y,0) := \frac{y}{r} q(r), \\
 \rho(x,y,0) &:= \gamma \left(1 - \frac{\gamma - 1}{2} q(r) \right)^{\frac{2}{\gamma -1}}, \quad \quad \\
 \pp(x,y,0) &:= \frac{\rho(x,y,0) \left(1 - \frac{\gamma - 1}{2} q(r) \right)^2 }{\gamma}. 
\end{alignat}
See Fig. \ref{fig:radexpwave} for a picture of initial and final density. 
As the flow remains smooth, at least for $\gamma = 3$, the entropy 
\begin{align}
 s := \ln \left(\frac{\pp}{\rho^{\gamma}} \right) 
\end{align}
remains constant. This constant is denoted by $s_0$. 

It is expected for this test case to monitor the $L^2$-norm of the average entropy error. We perform this exercise for two different (uniform) grids, one having 2048 elements, and the other having 8192 elements, and for quadratic and cubic ansatz functions. Entropy is monitored for all the adaptive DIRK methods available and compared against a BDF2 and BDF3 scheme, respectively, see Figs. \ref{fig:convg_radexpwave2048p2}-\ref{fig:convg_radexpwave8192p3}. Tolerance $\tol$ was chosen to be $\tol = 0.5$ for quadratics and 2048 elements, and $\tol = 10^{-1}$ for the other computations.
It can be clearly seen that the adaptive DIRK methods perform as well as the BDF schemes, except for the last test case, where only the Runge-Kutta method by Hairer and Wanner performs as good as BDF3. Nevertheless, the deviation of Cash's and Al-Rabeh's method from BDF3 is not too extreme. What can again be observed is the fact that the adaptive methods choose an appropriate time-step that is an order of magnitude larger than the one corresponding to the BDF schemes. (Note that the time-step size that corresponds to BDF is chosen in such a way that temporal resolution has minimal effect on the accuracy of the entropy. This is one of the requirements from the high order workshop.) Furthermore, for $t \rightarrow \infty$, the flow gets more and more trivial. This is, for all the methods, reflected in the time-step size that is increasing. 
The most expensive part of an implicit time integration method is the Newton steps. For this reason, we document the cumulated number of Newton iterations over time, including the rejected steps, see Figs. \ref{fig:convg_radexpwave2048p2}-\ref{fig:convg_radexpwave8192p3}. Even in this not so long-term run, it can be seen that there is 
a clear advantage of the 
 adaptive methods, especially for the methods by Hairer and Wanner and Al-Rabeh. Furthermore, the curves associated to the adaptive schemes have a much smaller slope than the curve associated to the BDF methods. For long-time runs, this constitutes a clear advantage.

\subsection{Navier-Stokes equations (Von K\'arm\'an vortex street)}

\begin{figure}
 \includegraphics[width=0.9\textwidth]{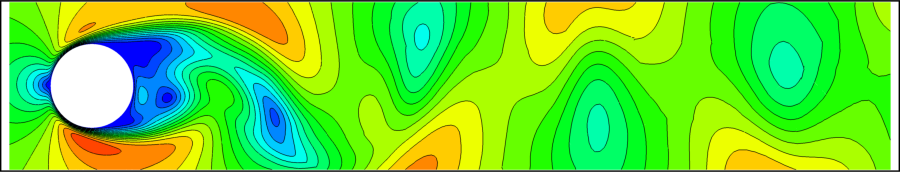} \\ \vspace{0.1cm} \\
 \includegraphics[width=0.9\textwidth]{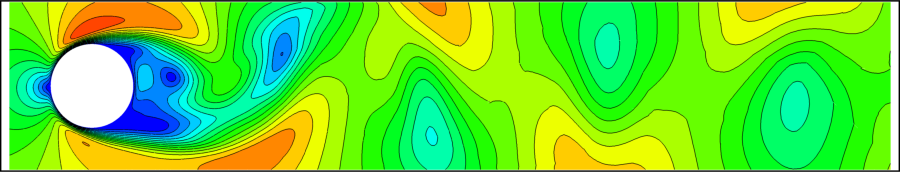}
 \caption{Von K\'arm\'an vortex street: Mach number at two instances.} \label{fig:vonkarman}
\end{figure}

\begin{figure}
 \begin{center}
  \includegraphics{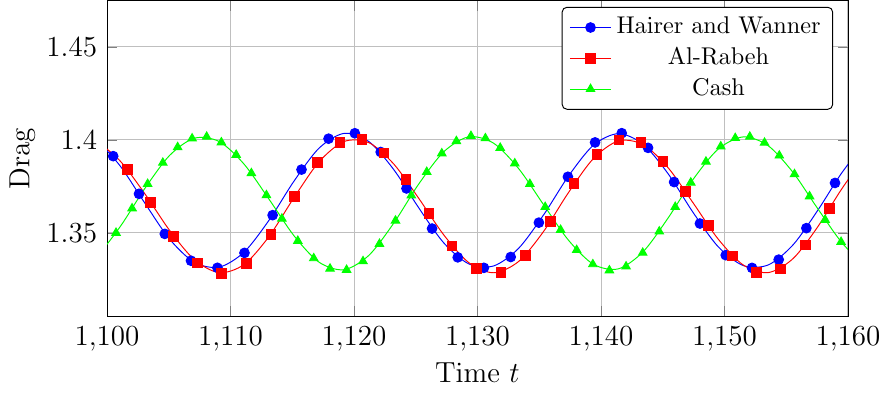}
  \caption{Von K\'arm\'an vortex street: Evolution of drag coefficient.}
  \label{fig:drag_vonkarman}
 \end{center}
\end{figure}

\begin{figure}
 \begin{center}
  \includegraphics{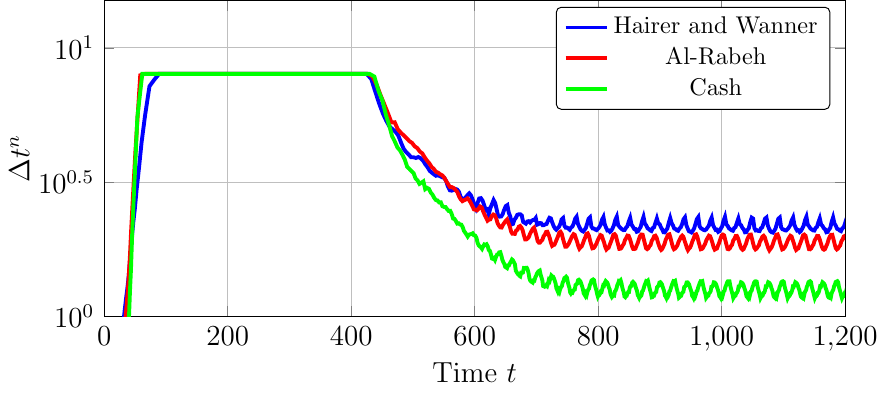}
  \caption{Von K\'arm\'an vortex street: Time-step evolution.}\label{fig:timestep_vonkarman}
 \end{center}
\end{figure}

The last numerical test case computes a von K\'arm\'an vortex street, see Fig. \ref{fig:vonkarman}. It is well-known that for $\Rey > 50$, flow around a circular cylinder gets unstable and the process of vortex shedding begins. 
We choose free stream Mach $\Ma$ and Reynolds number $\Rey$, respectively, as $\Ma = 0.2$ and $\Rey = 180$. Along the cylinder, we use no-slip boundary conditions, while in the farfield, we use a characteristic inflow/outflow boundary condition based on the freestream values. 
The employed mesh consists of 2916 elements and extends to 20 diameters away from the cylinder. We used the same mesh in our earlier work, see \cite{SchWoMay2012}. Computations are performed with cubic ansatz functions, and a tolerance of $\tol = 10^{-2}$. 
We choose a maximum time-step size of 8, i.e., we enforce $\Delta t^n \leq 8$, because otherwise, the physics of the flow are not correctly captured. Minimum time-step is chosen as $10^{-2}$, i.e., we enforce $\Delta t^n \geq 10^{-2}$. This prevents the time-step control from reaching too little values in the beginning of the flow. 
There is vast literature on this test case. For the free stream values as indicated one can, e.g., compute mean drag coefficients and Strouhal numbers. Reference values have been reported in literature, see Tbl. \ref{tbl:cyl_literature}. Next to these reference values, we have tabulated our computational values. Furthermore, in Fig. \ref{fig:drag_vonkarman}, we have plotted the drag evolution. Note that all the methods produce a periodic drag distribution with nearly the same period. However, as the unsteady process in this setting is introduced by the unstable nature of the flow, it depends on the startup phase when the process of vortex shedding begins. This explains why the drag coefficients of the different methods are shifted. 
In Fig. \ref{fig:timestep_vonkarman}, we have plotted the time-step evolution for the three different methods. One can see that in the beginning, the limiting of the time-step is really needed. It is only after initial instabilities have formed that the time-step is determined in a useful way. (Note that $\Delta t^n$ is periodic for all the methods, which resembles the periodic nature of the flow.) Again, it can be seen that Cash's method uses the smalles time-step, which is to be expected as it is the method with lowest nominal order. 


\begin{table}[!ht]
 \centering
 \begin{tabular}{|c|c|c|}
  \hline
  Experiment & $c_D$ & $\mathrm{Sr}$ \\ \hline
  Gopinath \cite{Gopinath2006} & 1.3406 & 0.1866 \\
  Henderson \cite{Henderson1995} & 1.336 & -\\
  Williamson \cite{Williamson96} & - & 0.1919 \\
  \hline
 \end{tabular}
 \begin{tabular}{|c|c|c|}
  \hline
  Time Discretization 	& $c_D$ & $\mathrm{Sr}$ \\ \hline
  Hairer and Wanner 	& 1.3666 & 0.1909 \\
  Al-Rabeh 		& 1.3653 & 0.1928\\
  Cash 			& 1.3672 & 0.1905\\
  \hline
 \end{tabular}
 \caption{Mean drag coefficients and Strouhal numbers from literature and computations}\label{tbl:cyl_literature}
\end{table}

%

\section{Conclusions and Outlook}\label{sec:conclusions}
In this work, we have developed a combination of a hybridized discontinuous Galerkin method and an embedded diagonally implicit Runge-Kutta method. We have shown numerical results that demonstrate how different Runge-Kutta methods perform. It seems that the method taken from the book by Hairer and Wanner \cite{HaiWan} performs best in most of our examples, whereas the method by Cash \cite{Cash} often needs a smaller time-step and seems to perform a somehow worse.

For the test problems considered in this paper, the combination of hybridized DG with embedded DIRK seems to perform very well. For problems with moving meshes, it has been recognized that space-time Galerkin methods work well \cite{WaPe13}. In the context of hybridized DG methods for incompressible flows, this has been investigated in \cite{RhebCo12,RhebCoVe13}. An interesting topic of future work is the extension of this to compressible flows. 
Future work will also include the investigation of multiderivative time integrators in the context of HDG as presented in \cite{Seal13}, because in principle, also the derivative is at hand in implicit methods. However, both stability and efficiency issues associated with this method have to be investigated. 
Furthermore, it has been recognized \cite{HaiWan, Bo07} that  in the context of singularly perturbed problems, e.g., for low Mach number flow, the convergence order of DIRK methods (and other Runge-Kutta methods including fully implicit ones, meaning the Butcher tableau is a dense matrix) can deteriorate. Fully implicit Runge-Kutta methods (e.g., Radau methods \cite{Bu64}) offer at least a partial remedy in that the deterioration is not that strong. Coupling this to the hybridized DG methods is left for future work. 

Adaptation in the temporal domain alone is obviously not enough to get the full potential of a method. Spatial adaptation is already available in the solver \cite{BaWoMa13}, and temporal adaptation has been investigated in this paper. Future work should therefore couple both ingredients in a sophisticated way to achieve maximum efficiency. One idea is to use an adjoint-based error indicator in both space and time to optimally design the spatial mesh. Adaptation in time can then be performed by a mixture of the adjoint and the DIRK time-step prediction.

\appendix
\section{Embedded DIRK schemes}\label{sec:appendix_dirk}
In this short appendix, we have listed the embedded DIRK schemes that we employ in our numerical results section:
The first tableau can be found in the classical book by Hairer and Wanner \cite{HaiWan}, its design order of accuracy is 4 and 3, respectively. 
\begin{equation}
     \begin{array}{c|cccccc}
         \frac{1}{4}    & \frac{1}{4}    & 0 & 0 & 0 & 0 & \vspace{0.02cm} \\ 
         \frac{3}{4}    & \frac{1}{2} & \frac{1}{4} & 0 & 0 & 0 & \vspace{0.02cm}\\
         \frac{11}{20}    & \frac{17}{50} & -\frac{1}{25} & \frac{1}{4} 
& 0 & 0 & \vspace{0.02cm}\\
         \frac{1}{2}& \frac{371}{1360} & -\frac{137}{2720} & 
\frac{15}{544} & \frac{1}{4} & 0 & \vspace{0.02cm}\\
         1            & \frac{25}{24} & -\frac{49}{48} & \frac{125}{16} 
& \frac{-85}{12} & \frac{1}{4} & \vspace{0.02cm}\\
     \hline
         \gamma_1& \frac{25}{24} & -\frac{49}{48} & \frac{125}{16} & 
-\frac{85}{12} & \frac{1}{4} & \vspace{0.02cm}\\
         \gamma_2& \frac{59}{48} & -\frac{17}{96} & \frac{225}{32} & 
-\frac{85}{12} & 0 & \vspace{0.02cm}\\
     \end{array}
     \label{eq:butcher:HW_45}
\end{equation}
The second tableau is due to Al-Rabeh \cite{Rabeh87}, its design order of accuracy is 4 and 3, respectively.
\begin{equation}
     \begin{array}{c|ccccc}
         0.4358665 &  0.4358665 &  0         &  0 & 0 & \\
         0.0323722 & -0.4034943 &  0.4358665 &  0 & 0 & \\
         0.9676278 & -0.3298751 &  0.8616364 &  0.4358665 & 0 & \\
         0.5641335 &  0.5575315 & -0.1930865 & -0.2361781 & 0.4358665 &  \\
     \hline
         \gamma_1      & 0.3153914 & 0.1846086 & 0.1846086 & 0.3153914 & \\
         \gamma_2& 0.6307827 & 0.1413538 & 0.2278634 & 0 & \\
     \end{array}
     \label{eq:butcher:AR_34}
\end{equation}
The last tableau we use is due to Cash \cite{Cash}, its design order of accuracy is 3 and 2, respectively. 
\begin{equation}
     \begin{array}{c|cccc}
         0.435866521508 &  0.435866521508 &  0         & 0         & \\
         0.717933260755 & 0.2820667320 &  0.435866521508 &  0         & \\
         1.0 & 1.208496649 &  -0.6443632015 & 0.435866521508 & \\
     \hline
         \gamma_1      & 1.208496649 &  -0.6443632015 & 0.435866521508 &  \\
         \gamma_2& 0.77263013745746 & 0.22736986254254 & 0.0 &  \\
     \end{array}
     \label{eq:butcher:Cash_33}
\end{equation}

\begin{figure}[h]
 \begin{center}
   \subfigure[Evolution of entropy error.]{\includegraphics{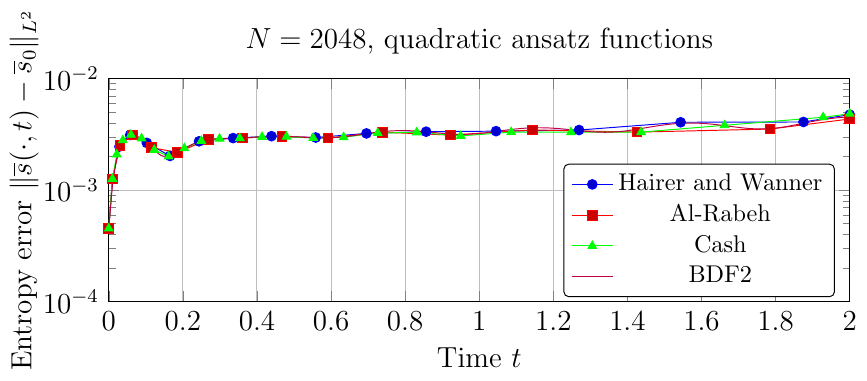}\label{fig:convg_radexpwave2048p2_ent}}
   \subfigure[Time-step evolution.]{\includegraphics{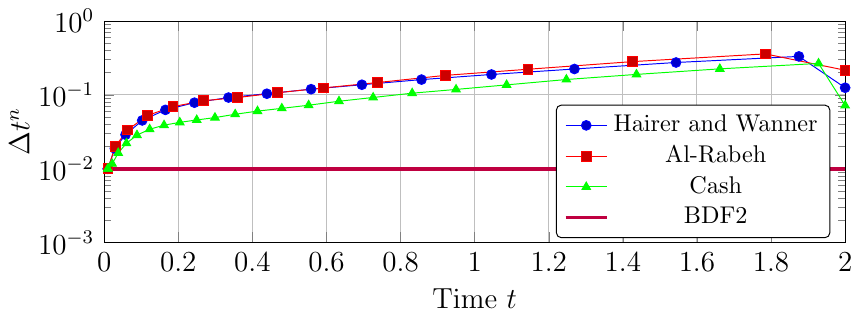}\label{fig:convg_radexpwave2048p2_dt}}
   \subfigure[Cumulated Newton steps.]{\includegraphics{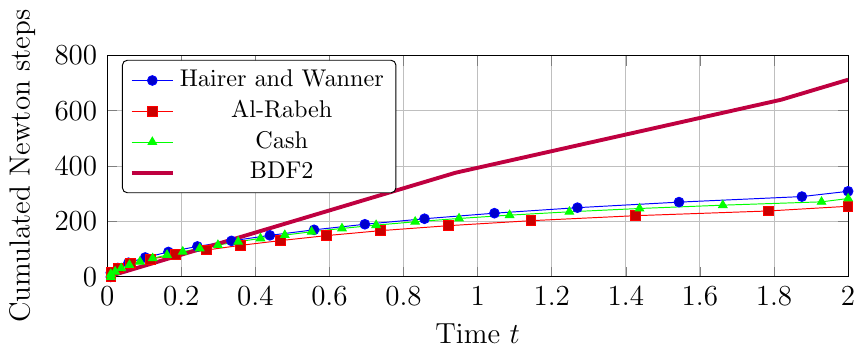}\label{fig:convg_radexpwave2048p2_newton}}
   \caption{Radial Expansion Wave: 2048 elements and quadratic ansatz functions.}\label{fig:convg_radexpwave2048p2}
 \end{center}
\end{figure}

\begin{figure}[h]
 \begin{center}
   \subfigure[Evolution of entropy error.]{\includegraphics{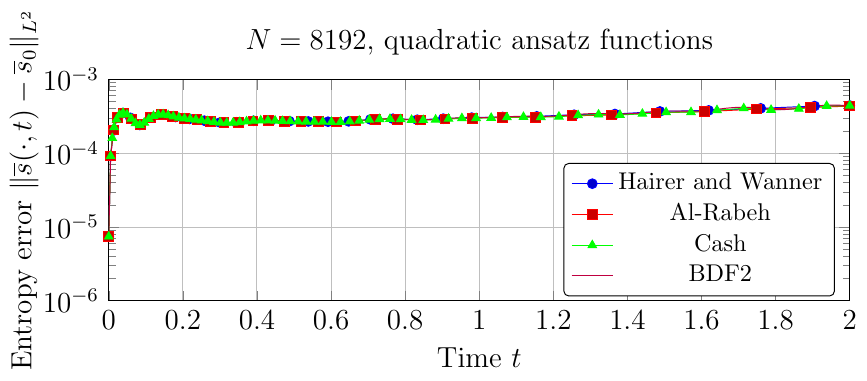}\label{fig:convg_radexpwave8192p2_ent}}
   \subfigure[Time-step evolution.]{\includegraphics{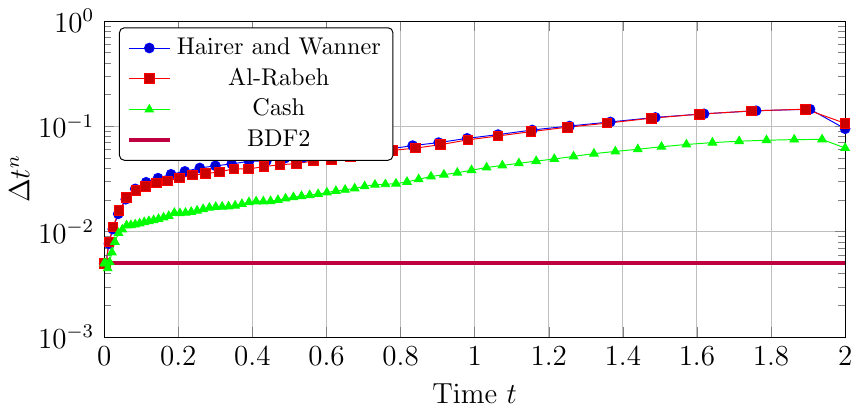}\label{fig:convg_radexpwave8192p2_dt}}
   \subfigure[Cumulated Newton steps.]{\includegraphics{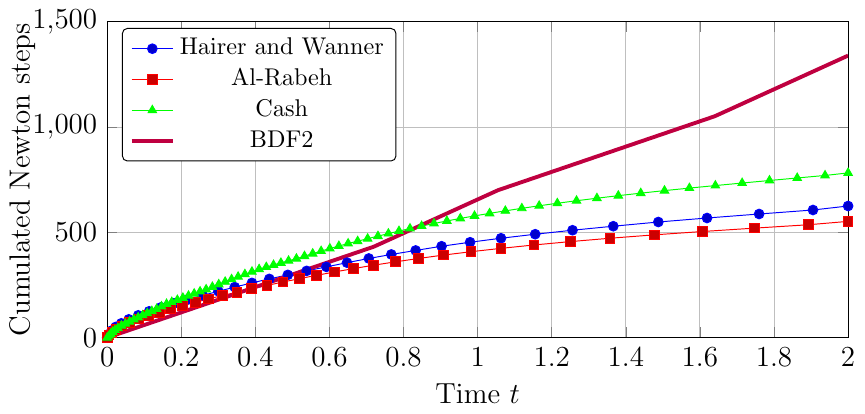}\label{fig:convg_radexpwave8192p2_newton}}
   \caption{Radial Expansion Wave: 8192 elements and quadratic ansatz functions.}\label{fig:convg_radexpwave8192p2}
 \end{center}
\end{figure}

\begin{figure}[h]
 \begin{center}
   \subfigure[Evolution of entropy error.]{\includegraphics{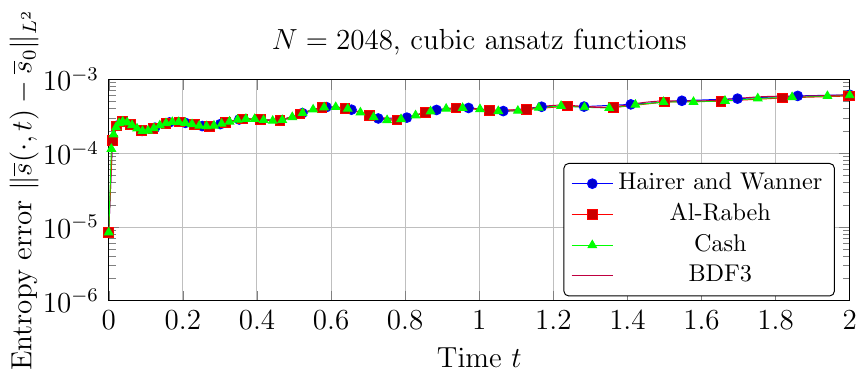}\label{fig:convg_radexpwave2048p3_ent}}
   \subfigure[Time-step evolution.]{\includegraphics{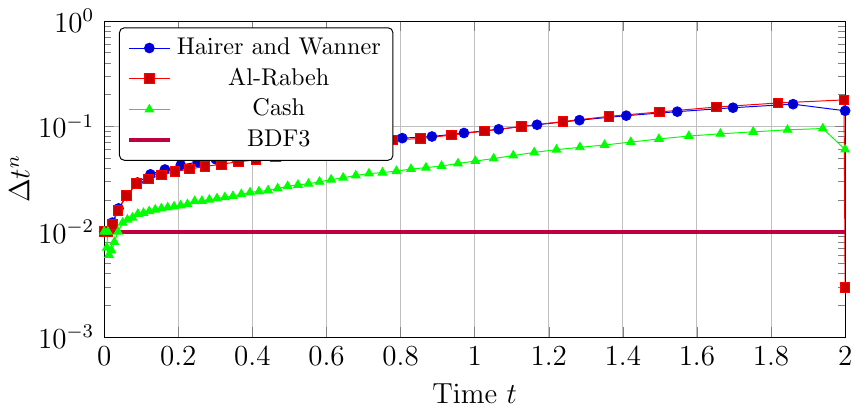}\label{fig:convg_radexpwave2048p3_dt}}
   \subfigure[Cumulated Newton steps.]{\includegraphics{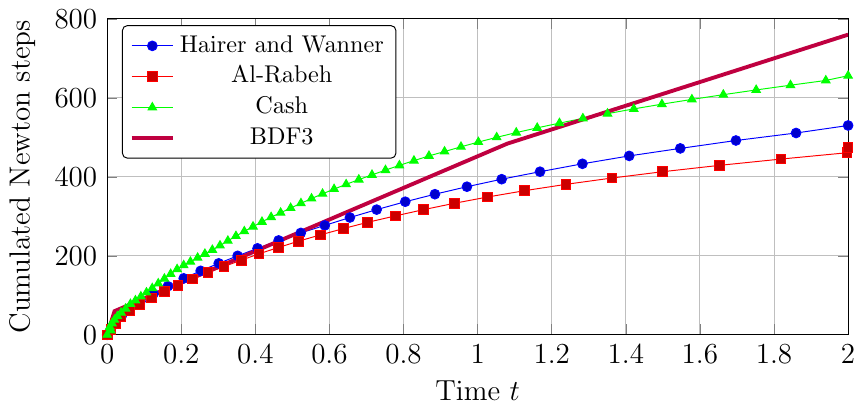}\label{fig:convg_radexpwave2048p3_newton}}
   \caption{Radial Expansion Wave: 2048 elements and cubic ansatz functions.}\label{fig:convg_radexpwave2048p3}
 \end{center}
\end{figure}

\begin{figure}[h]
 \begin{center}
   \subfigure[Evolution of entropy error.]{\includegraphics{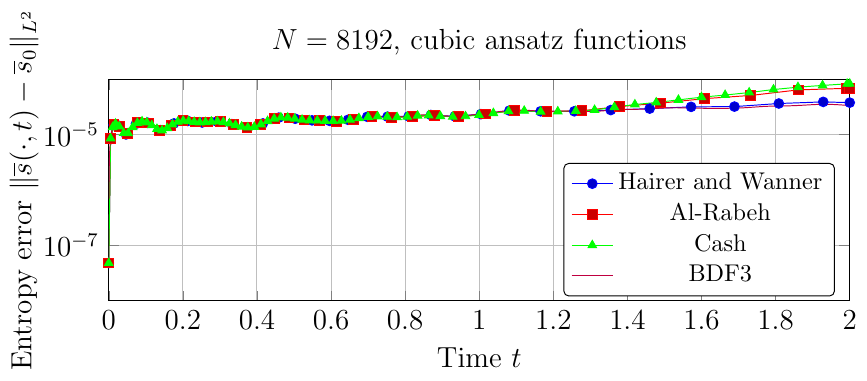}\label{fig:convg_radexpwave8192p3_ent}}
   \subfigure[Time-step evolution.]{\includegraphics{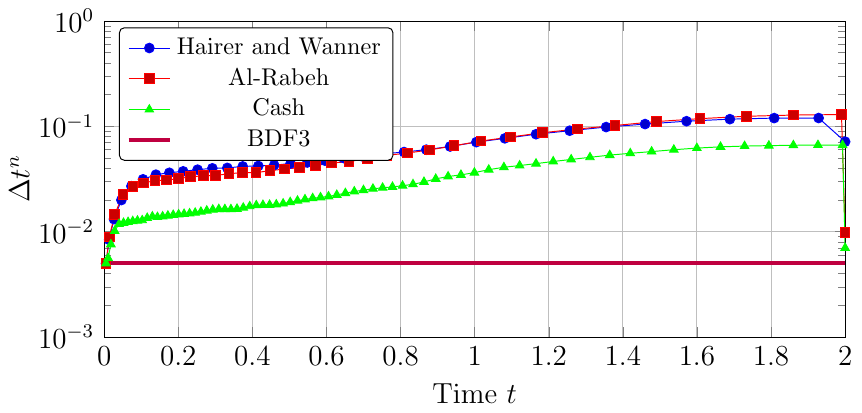}\label{fig:convg_radexpwave8192p3_dt}}
   \subfigure[Cumulated Newton steps.]{\includegraphics{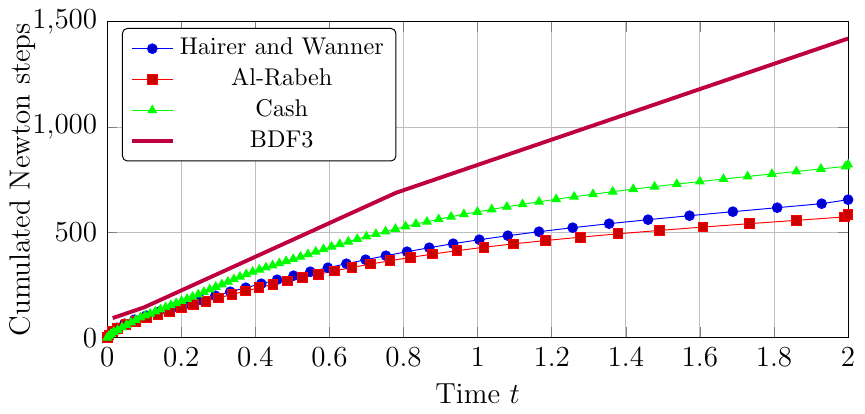}\label{fig:convg_radexpwave8192p3_newton}}
   \caption{Radial Expansion Wave: 8192 elements and cubic ansatz functions.}\label{fig:convg_radexpwave8192p3}
 \end{center}
\end{figure}

\bibliographystyle{elsarticle-num}
\bibliography{ListPaper}



\end{document}